\documentclass{amsart}[10pt]
\usepackage{amsmath}
\usepackage{amssymb} 
\usepackage{epsbox}
\usepackage{amsthm}
\usepackage{graphicx}

 \newtheorem{thm}{Theorem}
 \newtheorem{cor}{Corollary}

 \newtheorem{prop}{Propotition}
 \newtheorem{lem}{Lemma}
 \newtheorem{claim}{Claim}

 {\theoremstyle{definition}
 \newtheorem{rem}{Remark}}
 {\theoremstyle{definition}
  \newtheorem{defn}{Definition}}
 {\theoremstyle{definition}
 \newtheorem{exam}{Example}}

 \newcommand{\C}{ \mathcal{C}}
 \newcommand{\bk}{\mathbf{k}}

\begin{document}  

\title{Finite Thurston type orderings on dual braid monoids}
\author{Tetsuya Ito}

\begin{abstract} 
For a finite Thurston type ordering $<$ of the braid group $B_{n}$, we introduce a new normal form of a dual positive braid which we call the $\C$-normal form. This normal form extends Fromentin's rotating normal form and the author's $\C$-normal form of positive braids.
Using the $\C$-normal form, we give a combinatorial description of the restriction of the ordering $<$ to the dual braid monoids $B_{n}^{+*}$. We prove that the restriction to the dual positive monoid $(B_{n}^{+*},<)$ is a well-ordered set of order type $\omega^{\omega^{n-2}}$. 
\end{abstract}
\maketitle

\section{Introduction}

This is a subsequent of author's study of Thurston type orderings.
Thurston type orderings are left-invariant total orderings of the braid groups $B_{n}$, which is defined by using the hyperbolic geometry \cite{sw}. 

In the author's previous paper \cite{i}, the author provide a combinatorial description of a finite Thurston type ordering by introducing a new normal form of positive braids which we call the $\C$-normal form. Using the $\C$-normal form description, we proved that the restriction of a finite Thurston type ordering to the positive braid monoid $B_{n}^{+}$ is $\omega^{\omega^{n-2}}$.

There are other well-known submonoids of the braid groups, called the {\it dual braid monoids} denoted by $B_{n}^{+*}$. As in the positive braid monoids, the dual braid monoids also define a Garside structure of the braid groups and these two monoids are closely related each other \cite{b}.

The main theme of this paper is to study a finite Thurston type ordering of braid groups by using the dual braid monoid. Although the dual braid monoids and positive braid monoids have almost the same algebraic properties, sometimes the dual braid monoids provide much better result than the positive braid monoids. Fromentin studied the Dehornoy ordering, which is a special one of the finite Thurston type orderings, by the dual braid monoids in \cite{f1},\cite{f2},\cite{f3}. Using the dual braid monoid, he showed that every braid admits a $\sigma$-definitive quasi-geodesic, which has been conjectured for a long time. Thus it deserve to study the structure of finite Thurston type orderings on dual braid monoids.

We extend our code constructions in author's previous paper \cite{i} for the dual braid monoids. 
We present a combinatorial description of finite Thurston type orderings by introducing the new normal form of dual positive braids, which we call the $\C$-normal form.
The definition of the $\C$-normal form of dual positive braids is parallel to the definition of the $\C$-normal form of positive braids.
For a finite Thurston type ordering $<$, we define the code of a dual positive word, which is an sequence of parenthesized integers. The set of codes has two natural lexicographical orderings $<_{left}$ and $<_{right}$. The $\C$-normal form of a dual positive braid $\beta$ is defined as a word representative of $\beta$ which has the maximal code with respect to the ordering $<_{right}$. We denote the code of the $\C$-normal form of $\beta$ by $\C(\beta;<)$.
We also present an alternative definition of the $\C$-normal form, which is more algebraic and useful to actual computation.
Fromentin defined the normal form of dual positive braids called the {\it rotating normal form}. We introduce the {\it tail twisted rotating normal form} of dual positive braids, which extends the rotating normal form. We show that the $\C$-normal form is identified with the tail-twisted rotating normal form, thus the $\C$-normal form can also be regarded as an extension of the rotating normal forms. 

Our main result is the following.
\begin{thm}
\label{thm:main}
Let $<$ be a finite Thurston type ordering on $B_{n}$.
Then for $\alpha,\beta \in B_{n}^{+}$, $\alpha < \beta$ holds if and only if $\C(\alpha;<) <_{left} \C(\beta;<)$ holds.
\end{thm}

Using this $\C$-normal form description, we determine the order type of the restriction of a finite Thurston type ordering to the dual braid monoid $B_{n}^{+*}$.

\begin{thm}[Order types of finite Thurston type ordering]
\label{thm:ordertype}
Let $<$ be a finite Thurston type ordering on the braid group $B_{n}$.
Then the order type of $(B_{n}^{+*},<)$ is $\omega^{\omega^{n-2}}$.
\end{thm}
Finally we give a computational complexities. 
\begin{thm}
\label{thm:complexity}
Let $<$ be a finite Thurston type ordering of $B_{n}$. 
\begin{enumerate}
\item For each dual braid monoid word $W$ of length $l$, the $\C$-normal form of $W$ can be computed by the time $O(l^{2})$.
\item For a braid (not necessarily dual positive) $\beta \in B_{n}$, whether $\beta>1$ holds or not can be decided by the time $O(l^{2}n^{2})$ where $l$ is the word length of $\beta$ with respect to the band generators. 
\end{enumerate}
\end{thm}

Although the above results are almost parallel to the positive braid monoid case, but the proofs of theorems are essentially different.
 We need a new approach to prove the main theorem, because the proof we used in the positive braid monoid case does not work. 
In the positive braid monoid case, we did not compare two positive braids $W$ and $V$ directly. Instead, we constructed two new braids $\underline{W}<W$ and $V<\overline{V}$ which are easy to compare the ordering. We proved the $\C$-normal form description for positive braid monoid $B_{n}^{+}$ by comparing $\underline{W}$ with $\overline{V}$ instead of comparing $W$ with $V$. 
However, in the dual braid monoid case, the same approach fails. It is too complex to construct a word which corresponds to $\underline{W}$ although it might be possible for small $n$. This is because the dual braid monoids have the much larger number of generators so there are very many patterns of subwords which play the similar role of $\underline{W}$.

To overcome this combinatorial difficulty, we develop a geometric theory of $\C$-normal form. We show that the $\C$-normal form is a word which acts of the arc of a curve diagram most effectively (Proposition \ref{prop:keyprop}). This special property of the $\C$-normal form enable us to compute the initial segment of the image of the curve diagram and to compare the ordering.\\

\textbf{Acknowledgments.}
   The author would like to express his gratitude to his advisor Toshitake Kohno for his encouragement. He also acknowledge Patrick Dehornoy and Jean Fromentin for driving my attention to the study of the dual braid monoids and for stimulating conversation. 

\setcounter{thm}{0}

\section{Preliminaries}

In this section we briefly review the definition of dual braid monoids, finite Thurston type ordering and curve diagrams.
 
\subsection{Braid groups, dual braid monoids}
The braid group $B_{n}$ is defined by the presentation 
\[
B_{n} = 
\left\langle
\sigma_{1},\sigma_{2},\cdots ,\sigma_{n-1}
\left|
\begin{array}{ll}
\sigma_{i}\sigma_{j}=\sigma_{j}\sigma_{i} & |i-j|\geq 2 \\
\sigma_{i}\sigma_{j}\sigma_{i}=\sigma_{j}\sigma_{i}\sigma_{j} & |i-j|=1 \\
\end{array}
\right.
\right\rangle
.
\]
 We call the monoid generated by the positive generators $\{ \sigma_{1},\sigma_{2},\cdots,\sigma_{n-1} \}$ the {\it positive braid monoid} and denote it by $B_{n}^{+}$. An element of $B_{n}^{+}$ is called a {\it positive braid}. 
 
For $1\leq i<j \leq n$, let $a_{i,j}$ be a $n$-braid defined by
\[ 
a_{i,j}= (\sigma_{j-1}\sigma_{j-2}\cdots \sigma_{i+1})\sigma_{i} (\sigma_{i+1}^{-1} \cdots \sigma_{j-2}^{-1}\sigma_{j-1}^{-1}). 
\]
The braids $\{a_{i,j}\}$ are called the {\it band generators}.
The braid group $B_{n}$ has the following presentation using the band generators.
\[ B_{n} = 
\left\langle
a_{i,j} \
;
\left|
\begin{array}{l}
a_{i,j}a_{r,s}= a_{r,s}a_{i,j}, \;\; (j-s)(j-r)(i-s)(i-r)>0  \\  
a_{i,j}a_{i,k}= a_{j,k}a_{i,j} = a_{i,k}a_{j,k}, \;\;\; i<j<k
\end{array}
\right.
\right\rangle
.
\]

{\it The dual braid monoid} $B_{n}^{+*}$ is a monoid generated by positive band generators. We call an element of $B_{n}^{+*}$ a {\it dual positive braid}.

\begin{figure}[htbp]
 \begin{center}
\includegraphics[width=80mm]{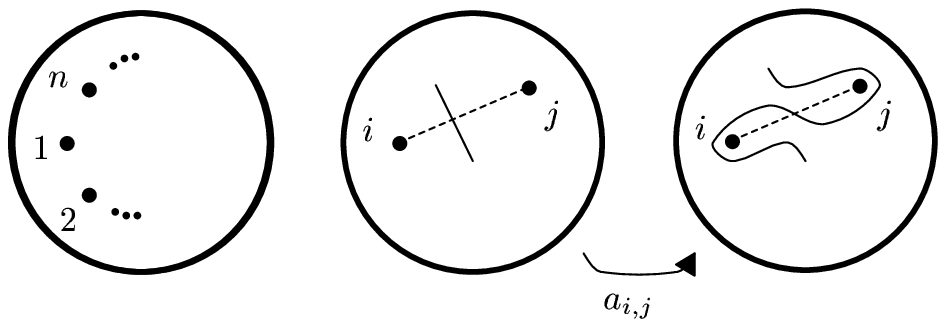}
 \end{center}
 \caption{Punctured disc $D_{n}$ and action of $a_{i,j}$}
 \label{fig:disc}
\end{figure}

In this paper, we regard the braid group $B_{n}$ as the mapping class group of the $n$-punctured disc $D_{n}$, which is the group of isotopy classes of homeomorphisms of $D_{n}$ whose restriction to the boundary are identity maps.
   To treat band generators, it is convenient to configure puncture points on the circle $x^{2}+y^{2}=\frac{1}{2}$, as shown in the figure \ref{fig:disc}. Let $e_{i,j}$ be an arc in the subdisc $x^{2}+y^{2} \leq \frac{1}{2}$, connecting the $i$-th and the $j$-th puncture points. Then the band generator $a_{i,j}$ corresponds to the isotopy class of the positive half Dehn twist along the arc $e_{i,j}$(See figure \ref{fig:disc}). 
   
\subsection{Finite Thurston type ordering and curve diagrams}

In this subsection we briefly review the definition of finite Thurston type orderings. 
We do not describe whole constructions of Thurston type orderings which uses the hyperbolic geometry \cite{sw}. Instead, we define finite Thurston type orderings by using a curve diagram. 
 The curve diagram approach of orderings first appeared in \cite{fgrrw} to give a geometric description of the Dehornoy ordering. The definition of curve diagrams and orderings given here is a natural generalization of \cite{fgrrw}, and carried out in \cite{rw} more general settings. 
      
A {\it curve diagram} $\Gamma$ is a union of oriented arcs $\Gamma_{1},\Gamma_{2},\cdots,\Gamma_{n-1}$ which satisfy the following properties.
\begin{enumerate}
\item Interiors of $\Gamma_{i}$ are disjoint to each other.
\item The initial point of $\Gamma_{i}$ lies in $(\bigcup _{j=1}^{i-1}\Gamma_{j}) \cup (\partial D_{n})$ and the end point of $\Gamma_{i}$ lies in (puncture points)$ \cup (\textrm{ interior of } \Gamma_{i})\cup (\bigcup _{j=1}^{i-1}\Gamma_{j}) \cup (\partial D_{n})$.
\item Each component of $ D_{n}\backslash \Gamma$ is a disc or an one-punctured disc.
\end{enumerate}

 We regard isotopic two curve diagrams as the same curve diagram.  
We say curve diagrams $\Gamma$ and $\Gamma'$ are tight if they form no bigons. A bigon is an embedded disc in $D_{n}$ whose boundary consists of two subarcs $\gamma \subset \Gamma_{i}$ and $\gamma' \subset \Gamma'_{j}$. If we equip a complete hyperbolic metric on $D_{n}$ and realize both $\Gamma$ and $\Gamma'$ as a union of geodesics, then $\Gamma$ and $\Gamma'$ are tight. Thus, we can always choose two curve diagrams so that they are tight.

For two distinct $\alpha,\beta \in B_{n}$, if we put the images $\alpha(\Gamma)$ and $\beta(\Gamma)$ are tight, then $\alpha(\Gamma)$ and $\beta(\Gamma)$ must diverge at some point. 
We say $\alpha <_{\Gamma} \beta$ holds if and only if $\beta(\Gamma)$ moves the left side of $\alpha(\Gamma)$ at the first divergence point.
   The relation $<_{\Gamma}$ defines a left-invariant total ordering on $B_{n}$. We call such a ordering a {\it finite Thurston type ordering}.
  
\subsection{Cutting sequence presentation of arcs}

  In this subsection we introduce a cutting sequence presentation of an embedded arc, which encodes a configuration of an arc into the sequence of signed integers. 

\begin{figure}[htbp]
 \begin{center}
\includegraphics[width=70mm]{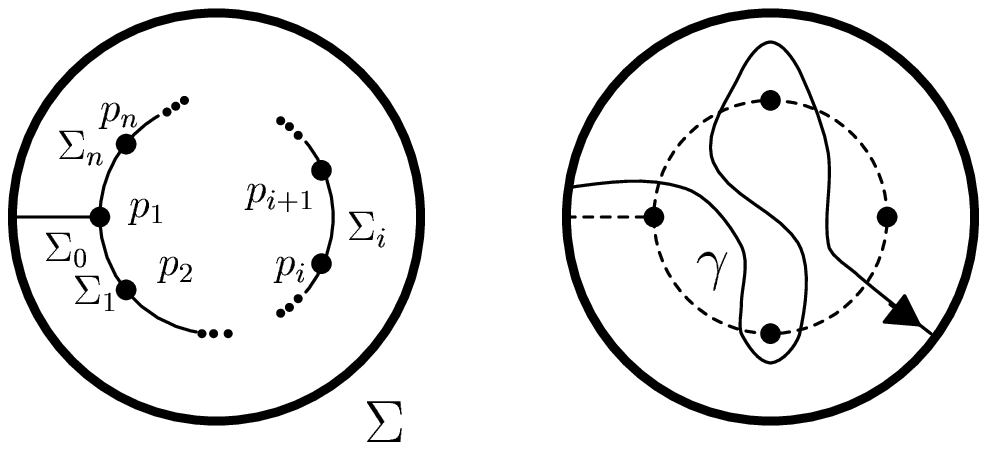}
 \end{center}
 \caption{Curve diagram $\Sigma$ and examples of cutting sequence}
 \label{fig:cut}
\end{figure}

Let $\Sigma_{0}$ be the straight arc which connects a point on the boundary and $p_{1}$, and define $\Sigma_{i} = e_{i,i+1}$ for $i=1,\cdots,n$, which are oriented from $p_{i}$ to $p_{i+1}$. Here we consider $e_{n,n+1}$ as $e_{n,1}$.
Let $\Sigma^{(n)}=\Sigma= \bigcup_{i=0}^{n} \Sigma_{i}$ be a curve diagram as shown in the figure \ref{fig:cut}.
We denote the disc bounded by the edge-path $\Sigma_{1}\cup \Sigma_{2}\cup \cdots \cup \Sigma_{n}$ by $D'$ or $D'_{n}$.

Let $\gamma$ be a properly embedded oriented arc in $D_{n}$. 
Isotope $\gamma$ so that $\gamma$ and the curve diagram $\Sigma$ transverse. Let $\{ p_{1},p_{2},\cdots,p_{k}\}$ be the intersection points of $\gamma$ with $\Sigma$. We define $C_{i}= + j$ (resp. $-j$) if $p_{i}$ lies on $\Sigma_{j}$ and the sign of the intersection at the point $p_{i}$ is positive (resp. negative). We call the sequence of signed integers $\{C_{1},C_{2},\cdots ,C_{k}\}$ {\it the cutting sequence presentation} of the arc $\gamma$. 
Since the complement $D_{n}\backslash \Sigma$ is a disc, up to isotopy the arc $\gamma$ can be reconstructed from its cutting sequence presentation.
We say a cutting sequence presentation is {\it tight} if the sequence contains no subsequences of the form $(\pm j, \mp j)$. Equivalently, a cutting sequence presentation is tight if and only if the corresponding arc is tight to $\Sigma$. 

\begin{exam}
The cutting sequence presentation of an arc $\gamma$ in right diagram of figure \ref{fig:cut} is $(+4,-1,+2,-4,+3,-2)$.
\end{exam}

\subsection{Normal orderings and examples}  
 
  In this subsection we introduce {\it normal finite Thurston type orderings}.
  Let $\Gamma = \bigcup_{i=1}^{n-1} \Gamma_{i}$ be a curve diagram. We say $\Gamma$ is {\it normal} if each $\Gamma_{i}$ is a properly embedded arc in $D_{n}$ and the tight cutting sequence presentation of $\Gamma_{i}$ is a length two sequence $(+n,-k(i))$. From the definition of a curve diagram, $k(i) \neq k(j)$ if $i\neq j$. Thus a normal curve diagram is represented by the permutation of $n-1$ integers $\bk=\{k(1),k(2),\cdots,k(n-1)\}$. 
We say a finite Thurston type ordering is {\it normal} if the ordering is defined by a normal curve diagram.

 Recall that two left invariant total orderings $<$ and $<'$ on a group $G$ are called {\it conjugate} if there exists an element $g \in G$ such that $f<h$ is equivalent to $fg<'hg$ for all $f,h \in G$. An element $g$ is called a {\it conjugating element} between $<$ and $<'$. In the braid group case, we can choose a conjugating element $g$ as a positive braid, hence as a dual positive braid. (\cite{i}, lemma 1). It is known that every finite Thurston type ordering is conjugate to a normal finite Thurston type ordering  \cite{sw}. Thus to obtain a description of a finite Thurston type ordering, it is enough to consider a normal finite Thurston type ordering.

We remark that this definition of normal finite Thurston type ordering is slightly different from that in our previous paper \cite{i}. We adopt this definition so that normal orderings are easy to treat in band generators. 
  
We close this section by giving examples of normal finite Thurston type orderings.

\begin{figure}[htbp]
\label{fig:circulardiagram}
 \begin{center}
\includegraphics[width=90mm]{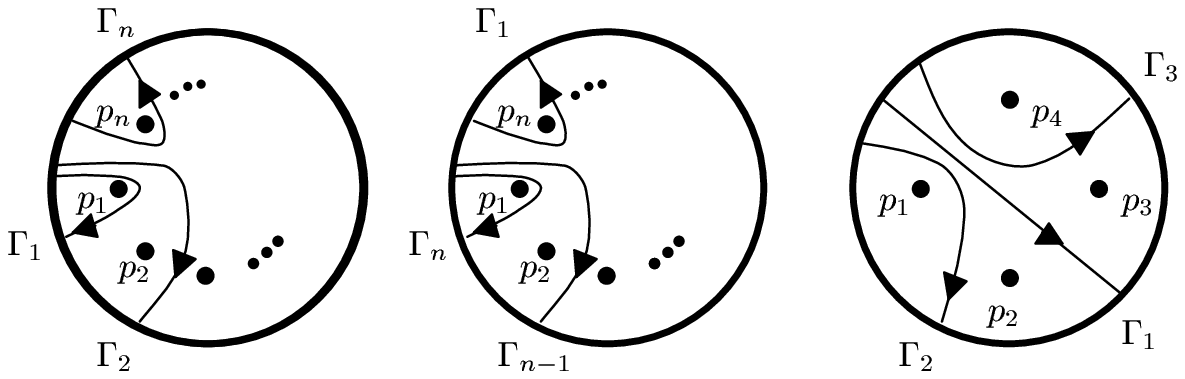}
 \end{center}
 \caption{Example of normal curve diagrams }
 \label{fig:example}
\end{figure}

\begin{exam}[Dehornoy ordering, Reverse of the Dehornoy ordering]
\label{exam:Dtype}
 Let $\Gamma$ be a curve diagram as shown in the right diagram in figure \ref{fig:example}. This curve diagram is represented by the trivial permutation $\bk=\{ 1,2,3,\cdots,n-1\}$.
The ordering defined by $\Gamma$ is called the {\it Dehornoy ordering} and denoted by $<_{D}$. Algebraically, the Dehornoy ordering $<_{D}$ is defined by $\alpha <_{D} \beta$ if and only if $\alpha^{-1}\beta$ admit a word representative which does not contain $\sigma_{1}^{ \pm 1},\cdots \sigma_{i-1}^{\pm1}, \sigma_{i}^{-1}$ but contains at least one $\sigma_{i}$.

Let $\Gamma'$ be a curve diagram the middle diagram of figure \ref{fig:example}, which corresponds to the permutation $\bk'=\{ n-1,n-2,\cdots,1\}$. The ordering defined by $\Gamma$ is called {\it the reverse of the Dehornoy ordering}. 

In some papers (especially, in Fromentin's paper \cite{f3} and Dehornoy's paper \cite{d2}), the Dehornoy ordering is used to represent the reverse of the Dehornoy ordering in this paper.
This is not a serious difference because the Dehornoy ordering and the reverse of the Dehornoy ordering are conjugate by the element $\Delta = (\sigma_{1}\sigma_{2}\cdots \sigma_{n-1})(\sigma_{1}\sigma_{2}\cdots \sigma_{n-2})\cdots(\sigma_{1}\sigma_{2})(\sigma_{1})$.  
\end{exam}

\begin{exam}
\label{exam:Ttype}
Let $\Gamma$ be a curve diagram in the right diagram of figure \ref{fig:example}, which is represented by the permutation $(2,1,3)$.
It is known that the ordering defined by $\Gamma$ is not conjugate to the Dehornoy ordering $<_{D}$. 
\end{exam}

\section{Code, $\C$-normal forms, tail-twisted $\Phi$-normal form}
In this section we provide algebraic and combinatorial aspects of our $\C$-normal form construction.

Throughout the paper, we use the following notational conventions.
For an integer $i$, we denote the unique integer $j$ which satisfies $j \equiv i \; (mod \; n)$, $1\leq j \leq n$ by $[i]$. Similarly, we use the notation $i^{*}$ to represent $n-i$ or $[n-i]=[-i]$. If we will not make a confusion, for example integers $i,j$ in a band generator $a_{i,j}$, we simply denote $i$ instead of the symbol $[i]$. 

For a triple of integers $n_{1},n_{2},n_{3}$, we denote $n_{1} \leq n_{2} \leq  n_{3}$ if the three puncture points $p_{[n_{1}]}, p_{[n_{2}]}$ and $p_{[n_{3}]}$ in $D_{n}$ are aligned anti-clockwise. For such triple $n_{1},n_{2}$ and $n_{3}$, the relation
\[ a_{n_{2},n_{3}}a_{n_{1},n_{2}} = a_{n_{1},n_{3}}a_{n_{2},n_{3}} = a_{n_{1},n_{2}}a_{n_{1},n_{3}} \]
holds. 
For a subset $I \subset \{1,2,\cdots,n\}$, we denote the submonoid of $B_{n}^{+*}$ generated by $\{a_{i,j}\: | \: i,j \in I \}$ by $B_{I}^{+*}$.
We will simply denote $B_{\{1,2,\cdots ,\widehat{i},\cdots,n\} }$ by $B_{(i)}^{+*}$.

Let $\phi$ be an inner automorphism of $B_{n}$ defined by $\phi(a_{i,j}) =a_{i+1,j+1}$. Then $B_{(i)}^{+*}$ is naturally identified with $B_{n-1}^{+*}$ by the map $\phi^{n-i}: B_{(i)}^{+*} \rightarrow B_{n-1}^{+*}$.
 
\subsection{Code and the $\C$-normal form}

   Let $<$ be a normal finite Thurston type ordering of $B_{n}$ and $\Gamma$ be a normal curve diagram which defines $<$, represented by the permutation $\bk=\{k(1),k(2),\cdots,k(n-1)\}$. Let $k=k(1)$.

For a dual braid positive word $W$, we define the code of $W$ with respect to the normal finite Thurston type ordering $<$, denoted by $\C(W;<)$ or simply $\C(W)$ by the following inductive way.

First we define the code of a dual positive $2$-braid word $W$. There is only one normal finite Thurston type ordering on $B_{2}$, namely, the Dehornoy ordering $<_{D}$. For a word $W=a_{1,2}^{p}$, we define the code of $W$ with respect to the Dehornoy ordering $<_{D}$ by the integer $\C(W;<_{D})=(p)$.

Assume that we have already defined the code of a dual positive $i$-braid word with respect to an arbitrary normal finite Thurston type ordering $<$ on $B_{i}$, for all $i<n$.

Let $W$ be a dual positive $n$-braid word. Take a factorization of $W$ into subwords
\[ W= A_{m}A_{m-1}\cdots A_{1}A_{0}A_{-1} \]
 where
\[
 A_{i} \in B_{(m^{*})}^{+*} \textrm{ for } i\geq 0, \; A_{-1} \in B_{\{1,2,\cdots,k\}}^{+*}\times B_{\{k+1,k+2,\cdots,n\}}^{+*}.
\] 
  Then the code $\C(W;<)$ (with respect to the above factorization) is defined as a sequence of codes with fewer strands
\begin{eqnarray*}
 \C(W;<) & = & ( \; \C( \phi^{m}(A_{m}); <_{D}), \cdots , \C ( \phi^{i}(A_{i});<_{D}), \\
& & \hspace{1cm} \cdots, \C(\phi^{1}(A_{1});<_{D}), \C( A_{0}; <_{res}), \C(A_{-1}; <) \;) 
\end{eqnarray*}
  where the last two terms $\C(A_{-1};<)$ and $\C(A_{0};<_{res})$ are defined by the following manner.
  
 Let $<_{res}$ be the restriction of the ordering $<$ to $B_{(n)}^{+*} = B_{n-1}^{+*}$. This is also a normal finite Thurston type ordering of $B_{n-1}$.
We define $\C( A_{0} ;<_{res})$ as the code of the $n-1$ dual braid monoid word $A_{0}$ with respect to the ordering $<_{res}$.
 
 Next we define the last term $\C(A_{-1};<)$. For $1 \leq j \leq n-2$, let $m_{j},M_{j}$ be integers defined by
\[
\left\{
\begin{array}{l}
m_{j} = \max \{1,k(i)+1 \:|\: i=1,2,\cdots,j, \; k(i)<k(j+1)\} \\
M_{j}= \min \{n,k(i) \:|\: i=1,2,\cdots,j,\; k(i)>k(j+1)\}.
\end{array}
\right.
\]
and $ I_{j} = \{  m_{j},m_{j+1},\cdots,M_{j} \}$.  
As in the word $W$, we decompose the word $A_{-1}$ as a product of subwords
\[ A_{-1} = X_{1}X_{2}\cdots X_{n-2}  \;\;\;\;\; (X_{j} \in B_{I_{j}}^{+*}) \]

For general word $A_{-1}$, such a form of decomposition might be impossible. We choose the first decomposition $W=A_{m}\cdots A_{-1}$ so that such a decomposition is possible. 

  Let $<_{j}$ be the restriction of the ordering $<$ on $B_{I_{j}}^{+*}$. By identifying $B_{I_{j}}^{+*}$ with $B_{M_{j}-m_{j}+1}^{+*}$ using the map $\phi^{-m_{j}+1}$, the ordering $<_{j}$ can be seen as a normal finite Thurston type ordering.  We have already defined the code of $X_{j}$ with respect to the ordering $<_{j}$. The last term $\C(A_{-1};<)$ is defined as a sequence of codes
  \[ \C(A_{-1};<)= ( \;\C(X_{1};<_{1}),\C(X_{2};<_{2}),\cdots,\C(X_{n-2};<_{n-2}) \;).\]
This completes the definition of codes.\\

  To define the code, it is not sufficient to choose a dual positive monoid word. we need to indicate the way to decompose $W$ into subwords. To indicate a subword decomposition of $W$, we sometimes use the symbol $|$.
For example, to express the decomposition of $W= a_{1,2}^{p+q}$ with $A_{-1}= a_{1,2}^{q}$ and $A_{0}= a_{1,2}^{p}$, we use a notation $W=a_{1,2}^{p}|a_{1,2}^{q}$.    
  
Let $Code(n*;<)$ be a set of all codes of dual $n$-braid monoid words with respect to the normal finite Thurston type ordering $<$. We define the two natural lexicographical orderings $<_{left}$ and $<_{right} $ on $Code(n*;<)$. To this end, we regard the code of a dual positive $n$-braid word as a left-infinite sequences of the code of subwords $(\cdots,a_{i},a_{i-1},\cdots,a_{1})$
, with only finite $a_{i}$ is non-trivial. This means, for example, we regard a code $(2)$ of a dual 2-braid word $a_{1,2}^{2}$ as a left infinite sequence of integers $(\cdots,0,0,2)$ not as a single integer $(2)$.

\begin{defn}[Lexicographical ordering $<_{left}$ and $<_{right}$]
Let $<$ be a normal finite Thurston type ordering.
The {\it Lexicographical ordering from left} $<_{left}$ and the {\it Lexicographical ordering from right} $<_{right}$ are total orderings on $Code(n*,<)$ define by the following inductive way.

For two codes $(a)$ and $(b)$ in $Code(2*,<_{D})$, we define $(a) <_{left} (b)$ if and only if $a<b$, and $(a) <_{right} (b)$ if and only if $a<b$.

Assume that we have already defined $<_{left}$ and $<_{right}$ for any normal finite Thurston type ordering $<$ on $B_{i}$ for all $i<n$. 
 
Let $\C = (\cdots,\C(A_{m}),\cdots, \C(A_{-1}) )$ and $\C' = (\cdots,\C(A'_{m}), \cdots,\C(A'_{-1})) $ be the codes of positive $n$-braid words $W$ and $W'$ with respect to the normal finite Thurston type ordering $<$. From definition of the codes, $\C(A_{-1}))$ and $\C(A'_{-1})$ are written by
\[
\left\{
\begin{array}{l} 
\C(A_{-1}) = ( \C ( X_{1}), \C (X_{2}),\cdots, \C(X_{n-2}) ) \\
\C(A'_{-1}) = ( \C ( X'_{1}), \C (X'_{2}),\cdots, \C(X'_{n-2}) ) 
\end{array}
\right.
\]
respectively. 

First we define $\C(A_{-1}) <_{left} \C(A'_{-1};<)$ if 
\[ \C(X_{j})=\C(X'_{j}) \textrm{ for } j<l  \textrm{  and } \C(X_{l}) <_{left} \C(X'_{l}) \textrm{ for some } l<n-2.\] 
Similarly, we define $ \C(A_{-1}) <_{right} \C(A'_{-1}; <)$ if 
\[ \C(X_{j})=\C(X'_{j}) \textrm{ for } j>l \textrm{ and } \C(X_{l}) <_{right} \C(X'_{l}) \textrm{ for some } l<n-2. \]

Now we define $\C <_{left} \C'$ if and only if 
\[
\C(A_{i})=\C(A'_{i}) \textrm{ for } i>l \textrm{ and } \C(A_{l}) <_{left} \C(A'_{l})\textrm{ for some } l \\
\]
and we define $\C <_{right} \C'$ if and only if 
\[
\C(A_{i})=\C_(A'_{i}) \textrm{ for } i<l \textrm{ and } \C(A_{l}) <_{right} \C(A'_{l}) \textrm{ for some } l.
\]

\end{defn}

Now we are ready to give a definition of the $\C$-normal form for normal finite Thurston type ordering.

\begin{defn}[The $\C$-normal form for a normal finite Thurston type ordering]
   Let $<$ be a normal finite Thurston type ordering on $B_{n}$.
For $\beta \in B_{n}^{+*}$, {\it the $\C$-normal form} of $\beta$ with respect to the ordering $<$ is the  dual braid monoid word representative of $\beta$ whose code is maximal among the set of codes of all dual monoid word representatives of $\beta$, with respect to the ordering $<_{right}$.
\end{defn}

   Since the number of dual braid monoid word representatives of $\beta$ is finite, the $\C$-normal form of a dual positive braid is always uniquely determined. For $\beta \in B_{n}^{+*}$, we denote the $\C$-normal form with respect to the ordering $<$ by $W(\beta;<)$ or simply $W(\beta)$. We denote its code by $\C(\beta;<)$ or $\C(\beta)$. 
   
 For a $\C$-normal form $W= A_{m}A_{m-1}\cdots A_{-1}$, we denote the integer $m$ by $a_{1}(W)$, and call it {\it the 1st address} of $W$. 

Let us denote the subword decomposition of $A_{m}$ by
\[ A_{m} = A^{(1)}_{m^{(1)}}A^{(1)}_{m^{(1)}-1}\cdots .\]
We denote the integer $m^{(1)}$ by $a_{2}(W)$ and call it {\it the $2$nd address} of $W$.
Similarly, using the decomposition of each $A^{(i-2)}_{m^{(i-2)}}$
\[ A^{(i-2)}_{m^{(i-2)}} = A^{(i-1)}_{m^{(i-1)}}A^{(i-1)}_{m^{(i-1)}-1} \cdots \]
we define the {\it $i$-th address} of $W$ by $a_{i}(W) =m^{(i-1)}$.

The depth of $\C$-normal form $d(W)$ is defined by
\[ d(W) = \min\{ i \geq 1 \:|\: A^{(i-1)}_{m^{(i-1)}} = a_{p,q}^{l} \} \]
The notion of address and depth will play an important role in our proof of theorem \ref{thm:main}.

Finally we extend the code and the $\C$-normal forms for arbitrary finite Thurston type orderings.
Let $<$ be a finite Thurston type ordering on $B_{n}$ which is conjugate to a normal ordering $<_{N}$, and $P$ be a dual positive word representative of a dual positive conjugating element between $<$ and $<_{N}$. 
\begin{defn}[The $\C$-normal form for a general finite Thurston type ordering]
Let $<$, $P$ and $<_{N}$ as the above.
For $\beta \in B_{n}^{+*}$, the $\C$-normal form of $\beta$ with respect to the ordering $<$ is defined by the dual positive word $N(\beta\cdot P;<_{N})P^{-1}$ where $N(\beta \cdot P; <_{N})$ is the $\C$-normal form of $\beta\cdot P$ with respect to the normal finite Thurston type ordering $<_{N}$. The code of $\beta$ with respect to the ordering $<$ is defined by $\C(\beta\cdot P;<_{N})$.
\end{defn}

From definition, the set of all codes $Code(n*,<)$ with respect to the ordering $<$ is a subset of $Code(n*,<_{N})$. The lexicographical orderings $<_{left}$ and $<_{right}$ on the set $Code(n*,<)$ are defined to be the restriction of the orderings $<_{left}$ and $<_{right}$ on the set $Code(n*,<_{N})$.

\subsection{Examples}
We provide examples of the code and the $\C$-normal forms.

\begin{exam}[The Dehornoy ordering]
\label{exam:dehornoy4}
First of all, we describe the $\C$-normal forms of the Dehornoy orderings.
Let us begin with the the $3$-braid case. Let $W$ be a dual positive $3$-braid word and take a decomposition of the word $W$ as
\[ W = \cdots A_{m}A_{m-1}\cdots A_{0}A_{-1}.\]
From definition, the above decomposition is explicitly written as  
\[ W = \cdots a_{1,2}^{p_{3}} a_{2,3}^{p_{2}} a_{1,3}^{p_{1}}a_{1,2}^{p_{0}}a_{2,3}^{p_{-1}}. \]
Thus, the code of $W$ is given by a sequence of integers
\[
\C(W;<_{D}) = (\cdots ,p_{3},p_{2},p_{1},p_{0},p_{-1}).
\]

For example, let us consider the word $W=a_{1,2}||a_{1,3}||$. Then its code is given by $\C(W) = (1,0,1,0,0)$. The word $W$ is not the $\C$-normal form. The $\C$-normal form of the braid $W$ is $W'=a_{1,3}||a_{2,3}$, whose code is $(1,0,1)$. The depth and the $1$st address of the $\C$-normal form $W'$are $d(W')=1$, $a_{1}(W') =1$, respectively.

Next we proceed to a more complicated example, the Dehornoy ordering of $B_{4}$.
To determine the code of a word $W$, we need to specify the way of decomposition
\[ W=A_{m}A_{m-1}\cdots A_{0}A_{-1}.\]

First observe that $m_{1}=2$, $M_{1}=4$, and $m_{2}=3$, $M_{2}=4$.
Thus, the restriction of $<$ to $B_{I_{1}}^{+*}$ and $B_{I_{2}}^{+*}$ are also the Dehornoy ordering, under the identification of $\phi^{3}$ and $\phi^{2}$ respectively.
Therefore the decomposition of $A_{-1}$ is written as 
\[ A_{-1} = X_{1}X_{2} = \cdots a_{3,4}^{q_{4}} a_{2,4}^{q_{3}} a_{2,3}^{q_{2}} a_{3,4}^{q_{1}} | a_{3,4}^{z} \] 
and its code is given by 
\[ \C(A_{-1};<) = ( \:(\cdots, q_{4},q_{3},q_{2},q_{1}), (z)). \]

 Next observe that the restriction of $<_{D}$ to $B_{\{1,2,3\}}^{+*}$ is also the Dehornoy ordering on $B_{3}$. So the word $A_{0}$ is explicitly written by
\[ A_{0} = \cdots  a_{2,3}^{p_{4}} a_{1,3}^{p_{3}} a_{1,2}^{p_{2}} a_{2,3}^{p_{1}}. \]

 Similarly, by definition, $\phi(A_{1})$ is the $\C$-normal form with respect to the Dehornoy ordering, so
$A_{1}$ is described as
\[ A_{1} = \cdots  a_{1,2}^{r_{4}} a_{2,4}^{r_{3}} a_{1,4}^{r_{2}} a_{1,2}^{r_{1}} \]
By the similar way, we can explicitly write each word $A_{i}$.

Now we proceed to consider a concrete example. Let $W=a_{3,4}||a_{1,2}a_{2,3}^{2}a_{1,2}|a_{3,4}$. That is, the decomposition is given by 
 \[ A_{2}=a_{3,4}, \; A_{1}= \varepsilon ,\; A_{0}=a_{1,2}a_{2,3}^{2}a_{1,2},\; A_{-1}=a_{3,4}. \]

For a sake of simplicity, we choose a subword decomposition of each $A_{i}$ so that they have the maximal code with respect to $<_{right}$ among all subword decomposition of $A_{i}$. 
Then, the code of each $A_{i}$ is given by
\[
\left\{
\begin{array}{l}
\C(\phi^{2}(A_{2}) ;<_{D}) = \C ( \phi^{2}(a_{3,4}|) ;<_{D}) = \C( a_{1,2}|; <_{D}) = (1,0) \\
\C(\phi^{1}(A_{1}); <_{D}) = \C ( \phi^{1}(\varepsilon); <_{D}) = \C( \varepsilon ; <_{D}) = (0) \\ 
\C( A_{0}); <_{res}) = \C( a_{1,2}|a_{2,3}^{2}||a_{1,2}|; <_{D}) =(1,2,0,1,0) \\
\C( A_{-1}; <) = (\C(\phi^{3}(\varepsilon); <_{1}=<_{D}),\C( \phi^{2}(a_{3,4});<_{2}=<_{D})) = ((0),(1)).\\
\end{array}
\right.
\]

Summarizing, we conclude that the code of $W$ is given by
\[ \C(W;<_{D}) = ((1,0),(0),(1,2,0,1,0),((0),(1))). \]

As is easily observed, $W$ is not the $\C$-normal form.
 The $\C$-normal form of the word $W$ is given by 
\[ W'= a_{1,2}a_{2,4}a_{1,4}||a_{2,4}a_{3,4}^{2} \]
and its code is
\[ \C(W';<_{D}) = (\: (1,1,1,0),(0),((1,0,0),(2))\: ). \]
The depth and addresses of $W'$ are given by
\[ d(W')=2, \; a_{1}(W')=1,\; a_{2}(W')=2. \]
\end{exam} 

\begin{rem}
\label{rem:difference}
Although the positive braid monoid $B_{n}^{+}$ is a submonoid of the dual braid monoid $B_{n}^{+*}$ and the definition of the $\C$-normal forms of $B_{n}^{+*}$ are very similar to that of $B_{n}^{+}$, the notion of $\C$-normal forms in positive and dual braid monoids are completely different. 
That is, in general the $\C$-normal form of a positive braid, which is defined in \cite{i}, are not the $\C$-normal form as a dual positive braid. Thus, the $\C$-normal form of a dual positive braid is not a simple extension of the $\C$-normal form of a positive braid. 

Let us explain this phenomenon by showing an example.
Let us consider the positive $4$-braid word $W = a_{3,4}a_{1,2}a_{2,3}^{2}a_{1,2}a_{3,4} = \sigma_{3}\sigma_{1}\sigma_{2}^{2}\sigma_{1}\sigma_{3}$. 
As we have shown in the example \ref{exam:dehornoy4}, $W$ is not the $\C$-normal form of as a dual positive braid. However, if we regard $W$ as a positive braid, then the word $W=\sigma_{3}\sigma_{1}\sigma_{2}^{2}\sigma_{1}\sigma_{3}$ is the $\C$-normal form with respect to the Dehornoy ordering (see \cite{i}).
\end{rem} 
   
\begin{exam}
\label{exam:thurston4}
Next we study the Thurston type ordering $<$ in example \ref{exam:Ttype}, which corresponds to the permutation $(2,1,3)$.

First of all, let us decompose a dual positive word $W$ as in the definition.
\[ W=A_{m}A_{m-1}\cdots A_{1}A_{0}A_{-1}.\]

We investigate each word $A_{i}$ and its code $\C(A_{i})$.
Observe that $I_{1}=\{(1,2),(2,1)\}$ and $I_{2}=\{(3,4),(4,3)\}$. Thus, the decomposition of the word $A_{-1}=X_{1}X_{2}$ is given by $A_{-1}= a_{1,2}^{p}a_{3,4}^{q}$, and the code of $A_{-1}$ is given by $\C(A_{-1};<) =(p,q)$.

Next we examine the code of $A_{0}$. Since the restriction of the ordering $<$ to $B_{\{1,2,3\}}^{+*}$ is the reverse of the Dehornoy ordering of $B_{3}^{+*}$ under the identification $\phi^{3}$, the word $A_{0}$ is written as 
\[ A_{0}= \cdots a_{1,2}^{p_{4}}a_{2,3}^{p_{3}}a_{1,3}^{p_{2}}a_{1,2}^{p_{1}}\]
and its code $\C(A_{0};<_{res})$ is given by 
\[ \C(A_{0};<_{res}) = (\cdots, p_{m},\cdots, p_{2},p_{1}).\]

The description of other words and codes are the same as the Dehornoy ordering case.

Now we give a concrete example. Let $W=a_{3,4}a_{1,2}a_{2,3}^{2}a_{1,2}a_{3,4} $ be a dual positive word appeared in the example \ref{exam:dehornoy4}. We choose the decomposition into subword by $W=a_{3,4}||a_{1,2}a_{2,3}^{2}|a_{1,2}a_{3,4}$. Then, the codes of the subwords $A_{i}$ are given by
\[
\left\{
\begin{array}{l}
\C(\phi^{2}(A_{2}); <_{D}) = \C(\phi^{2}(a_{3,4});<_{D})=\C ( a_{1,2}|;<_{D}) = (1,0) \\
\C(\phi^{1}(A_{1}); <_{D}) = \C ( \phi^{1}(\varepsilon ) ; <_{D}) = \C(\varepsilon;<_{D}) = (0) \\
\C( A_{0};<_{D}) = \C( a_{1,2}|a_{2,3}^{2}||;<_{res}) = (1,2,0,0) \\
\C( A_{-1};<_{D}) = \C( a_{1,2}|a_{3,4} ; <) = (1,1).
\end{array}
\right.
\]
 Summarizing, we conclude that the code of $W$ is given by 
 \[\C(W;<) = ((1,0),(0),(1,2,0,0),(1,1)).\]
This is not the $\C$-normal form.
The $\C$-normal form of $W$ with respect to the ordering $<$ and its code are given by  
\[ W' = a_{1,2}a_{2,4}^{2}||a_{1,2}a_{3,4}^{2},\; \C( W' ;<) = ( (1,2,0),(0),(1,2) ) \] 
respectively. The depth, and addresses of $W'$ is given by 
\[ d(W')= 2, \; a_{1}(W')= 1,\; a_{2}(W') = 1.\]
\end{exam}

\subsection{Tail twisted rotating normal form}

In this section we provide an alternative definition of the $\C$-normal form as an extension of Fromentin's rotating normal form, which is more algebraic and is very useful when we would like to compute.

The construction of the rotating normal form is based on the fact that $B_{n}^{+*}$ defines a Garside group structure \cite{bkl}. We do not give a detailed exposition of a Garside group structure, because we do not need the whole properties of a Garside structure. The property we use is the existence of the maximal right divisor. 

For $\beta \in B_{n}^{+*}$ and $I \subset \{ 1,2,\cdots, n \}$, the maximal right divisor of $\beta$ in $B_{I}^{+}$ is an element of $\beta' \in B_{I}^{+*}$ characterized by the following two properties.
\[
\left\{
\begin{array}{l}
\beta \beta'^{-1}\in B_{n}^{+*}. \\
\textrm{For } \beta'' \in B_{I}^{+*}, \textrm{ if } \beta \beta''^{-1} \in B_{n}^{+*}, \textrm{ then }\beta'\beta''^{-1} \in B_{I}^{+*}.
\end{array}
\right.
\]
We denote the maximal right divisor of $\beta$ in $B_{I}^{+*}$ by $\beta \wedge B_{I}^{+*}$.
The maximal right divisor always uniquely exists and is easily calculated. 

For a dual positive $n$-braid $\beta$, the rotating decomposition of $\beta$ is a factorization of $\beta$ \[ \beta = \beta_{m}\beta_{m-1}\cdots\beta_{0}\beta_{-1}\] 
where $\beta_{i}$ is defined by the following inductive formula.
\[
\left\{
\begin{array}{l}
\beta_{-1} = \beta \wedge B_{(1)}^{+*} \\ 
\beta_{i} = (\beta\cdot(\beta_{0}^{-1}\beta_{1}^{-1}\cdots \beta_{i-1}^{-1})) \wedge B_{(i^{*})}^{+*} \: : i\geq 0 . 
\end{array}
\right.
\]
 
Let us identify $B_{(i^{*})}^{+*}$ with $B_{n-1}^{+*}$ by the homomorphism $\phi^{i}$. Then we can iterate the rotating decomposition for each $\beta_{i}$ until each monoid is generated by one element. Thus,  iterating the rotating decomposition, we finally obtain the unique word representative $N(\beta)$ for each dual braid braid $\beta \in B_{n}^{+*}$.
We call this word the {\it rotating normal form} of $\beta$ and denote it by $N(\beta)$.

\begin{rem}
  The above description of the rotating normal form is slightly different from Fromentin's, because Fromentin used the reverse of Dehornoy ordering.  
\end{rem}

Now we generalize the above construction of the rotating normal form. 
Let $\bk =\{k(1),k(2),\cdots,k(n-1)\}$ be a permutation of $n-1$ integers $\{1,2,\cdots,n-1\}$. 
We define the $\bk$-tail twisted rotating normal form of a dual positive $n$-braid $\beta$, denoted by $\widetilde{N}(\beta;\bk)$, by the following inductive way. 

As in the definition of code, we begin with the $B_{2}^{+*}$ case. 
In this case, $\bk$ is the trivial permutation. We simply define the $\bk$-tail twisted rotating normal form of a dual positive $2$-braid $\beta=a_{1,2}^{p}$ by the word $\widetilde{N}(\beta;\bk) = a_{1,2}^{p}$.

Assume that we have already defined the $\bk$-tail twisted rotating normal form of a dual positive $i$-braid for all permutations $\bk \in S_{i-1}$ and $i<n$. Here $S_{i}$ denotes the degree $i$ symmetric group.
Then for $\bk \in S_{n-1}$ we define the $\bk$-tail twisted rotating normal form of a dual positive $n$-braid $\beta$ as follows.  

First we define the permutations $\bk_{0},\bk_{1},\cdots,\bk_{n-2}$. 
The permutation $\bk_{0}\in S_{n-2}$ is defined by
\[ \bk_{0} = \{ k(j_{1})-1,k(j_{2})-1,\cdots, k(j_{n-2})-1 \: | \: j_{i}< j_{i+1},\; k(j_{i}) > 1\}. \]

Let $m_{j}, M_{j}$ be integers and $I_{j} \subset \{1,2,\cdots,n\}$ as in the definition of codes (See 3.1).
The permutation $\bk_{j} \in S_{M_{j}-m_{j}}$ is defined by 
\[ \bk_{j} = \{ k(i_{1})-m_{j},\cdots, k(i_{M_{j}-m_{j}})-m_{j}\: |\: j<i_{q}< i_{q+1},\; m_{j} \leq k(i_{q}) \leq M_{j}\} .\]

Using these permutations, we consider the twisted version of the rotating decomposition. 
First of all, let us define $\beta_{Tail},\beta_{T_{0}}$ and $\beta_{Main}$ as  
\[
\left\{
\begin{array}{l}
\beta_{Tail} = \beta \wedge \left( B_{ \{ 1 , 2, \cdots, k(1) \} }^{+*} \times B_{\{ k_(1)+1,\cdots ,n \}}^{+*} \right) \\
\beta_{T_{0}}= (\beta \beta_{Tail}^{-1}) \wedge  B_{(n)}^{+*} \\
\beta_{Main}= \beta \beta_{Tail}^{-1}\beta_{T_{0}}^{-1}.
\end{array}
\right.
\]

We decompose $\beta_{Tail}$ as $\beta_{Tail} = \beta_{T_{1}}\beta_{T_{2}}\cdots\beta_{T_{n-2}}$ where $\beta_{T_{i}}$ is defined by

\[
\left\{
\begin{array}{lll}
\beta_{T_{n-2}} & = & \beta_{Tail} \wedge B_{I_{n-2}}^{+*} \\
\beta_{T_{n-j}} & = & (\beta_{Tail}\cdot \beta_{T_{n-2}}^{-1}\beta_{T_{n-3}}^{-1}\cdots \beta_{T_{n-j+1}}^{-1}) \wedge B_{I_{n-j}}^{+*}.
\end{array}
\right.
\]

Now each $B_{I_{j}}^{+*}$ is identified with $B_{M_{j}-m_{j}+1}^{+*}$ by $\phi^{-m_{j}+1}$. Using this identification, we denote the $\bk_{j}$-tail twisted rotating normal form of the braid $\beta_{T_{j}}$ by $\widetilde{N}(\beta_{T_{j}};\bk_{j})$. From an inductive assumption, we have already defined these tail twisted rotating normal forms. 
 The $\bk$-tail twisted rotating normal form of $\beta$ is defined by 
\[
\widetilde{N} (\beta;\bk) = N(\beta_{Main})\widetilde{N}(\beta_{T_{0}};\bk_{0})\cdots\widetilde{N}(\beta_{T_{n-2}};\bk_{n-2}).
\]
where we denote the rotating normal form of the braid $\beta_{Main}$ by $N(\beta_{Main})$. 

For the trivial permutation $\bk$, the $\bk$-tail twisted rotating normal form is merely Fromentin's usual rotating normal form. Thus, the $\bk$-tail twisted rotating normal form is an extension of the rotating normal form.

Now we show that the tail twisted rotating normal and the $\C$-normal form are the same.

\begin{prop}
\label{prop:twisted}
Let $<$ be a normal finite Thurston type ordering of $B_{n}$ represented by the permutation $\bk$.
Then the $\C$-normal form with respect to $<$ and the $\bk$-tail twisted rotating normal form coincide. 
\end{prop}
\begin{proof}

For a dual positive $n$-braid $\beta$, let $W=A_{m}A_{m-1}\cdots A_{-1}$ be the $\C$-normal form of $\beta$ and $A_{-1}=X_{1}X_{2}\cdots X_{n-2}$ be the decomposition of $A_{-1}$. Similarly, let $\beta_{Main},\beta_{Tail}$ and $\beta_{T_{j}}$ be dual positive braids as in the definition of the $\bk$-tail twisted rotating normal form and $\beta_{Main} =\beta_{l}\cdots\beta_{0}\beta_{-1}$ be the rotating decomposition of $\beta_{Main}$.
From the definition of $\beta_{T_{0}}$, both $\beta_{-1}$ and $\beta_{0}$ must be empty words. Since the subword $A_{-1}$ is chosen so that the code $\C(A_{-1})$ is maximal among the all subword decomposition of $\beta$, we obtain that $ A_{-1} = \beta \wedge \left( B_{ \{ 1 , 2, \cdots, k(1) \} }^{+*} \times B_{\{ k(1)+1,\cdots ,n \}}^{+*} \right) = \beta_{Tail}$.
By the similar arguments, we obtain $\beta_{T_{0}}=A_{0}$, $\beta_{i}=A_{i}$ and $X_{i}=\beta_{T_{i}}$. Thus, inductive argument give the desired result.
\end{proof}

\begin{exam}
Let $\beta$ be a dual positive 4-braid $\beta=a_{3,4}a_{1,2}a_{2,3}^{2}a_{1,2}a_{3,4}$, which appeared in example \ref{exam:thurston4} and $\bk=(2,1,3) \in S_{3}$. Then,
\[
\left\{
\begin{array}{l}
\beta_{Tail} = \beta \wedge B_{\{1,2\}}^{+*} \times B_{\{3,4\}}^{+*} = a_{1,2}a_{3,4}^{2} \\
\beta_{T_{0}} = (\beta\beta_{Tail}^{-1}) \wedge B_{\{1,2,3 \}}^{+*} = \varepsilon \\
\beta_{Main} = \beta_{1}=(\beta\beta_{Tail}^{-1}\beta_{T_{0}}^{-1}) \wedge B_{\{1,2,4\}}^{+*} = a_{1,2}a_{2,4}^{2}.
\end{array}
\right.
\]
Therefore $\beta_{T_{2}}=a_{1,2}a_{3,4}^{2} \wedge B_{\{3,4\}}^{+*} = a_{3,4}^{2}$ and $\beta_{T_{1}} = a_{1,2}$.
Now all permutations $\bk_{0}$, $\bk_{1}$ and $\bk_{2}$ are trivial, so all twisted rotating normal form for such permutations are usual rotating normal form. Thus,
\[
\left\{
\begin{array}{l}
\widetilde{N}(\beta_{T_{1}};\bk_{1}) = a_{1,2}, \; \widetilde{N}(\beta_{T_{2}};\bk_{2}) = a_{3,4}^{2}\\
\widetilde{N}(\beta_{T_{0}};\bk_{0}) = \varepsilon, \; N(\beta_{Main}) = a_{1,2}a_{2,4}^{2}.
\end{array}
\right.
\]
  
Thus we conclude that $\bk$-tail twisted rotating normal form of $\beta$ is
\[ \widetilde{N}(\beta;\bk)= (a_{1,2}a_{2,4}^{2})(\varepsilon)((a_{1,2})(a_{3,4}^{2})) = a_{1,2}a_{2,4}^{2}a_{1,2}a_{3,4}^{2}.\]
This is of course the same as the $\C$-normal form $\C(W;<)$ computed by example \ref{exam:thurston4}.
\end{exam}

\section{Proof of theorems and applications}

In this section we prove our main results and give applications.

\subsection{Principal part of arcs and $\C$-normal forms}

In this subsection, we introduce the principal part of a $\C$-normal form (or, a cutting sequences), which plays a central role of the proof of theorem \ref{thm:main}. 
Let $\Sigma^{(n)}=\Sigma= \bigcup_{i=0}^{n} \Sigma_{i}$ be a curve diagram in $D_{n}$ and $D' =D'_{n}$ be a subdisc of $D_{n}$ bounded by the edge-path $\Sigma_{1}\cup \cdots \Sigma_{n}$, defined in subsection 2.3. See figure \ref{fig:cut} again. For a dual positive braid $\beta$, we call the (tight) cutting sequence presentation of the arc $\beta(\Gamma_{1})$ simply the (tight) cutting sequence presentation of $\beta$. 

First we observe that for $\beta \in B_{n}^{+*}$, the tight cutting sequence presentation of $\beta$ is the following form.
\[
(+0,+0,\cdots, +0, \underbrace{+p,-q}_{\sharp},\cdots )
\]
 We call the subarc of $\beta(\Gamma_{1})$ corresponding to the subsequence $\sharp$ (or, the subsequence $\sharp$ itself) {\it the first principal part} of $\beta$. Geometrically, this is the initial subarc of $\beta(\Gamma_{1})$ which is cut by the disc $D' = D'_{n}$.  
 
We define the second, the third, $\cdots$, the $(n-2)$-th principal part of $\beta$ as follows.
Let $W=W(\beta;<)$ be the $\C$-normal form of $\beta$ and $a=a_{1}(W)$ be the $1$st address of $W$. We fill the $(n-a)$-th puncture point $p_{n-a}$ and shift a numbering of puncture points by $a$. That is, we choose a new numbering of the punctures defined by $p_{i} \rightarrow p_{i+a}$. After this operation, the obtained disc is naturally identified with the $n-1$ punctured disc $D_{n-1}$. 
Now, we consider the tight cutting sequence presentation of arc $\beta(\Gamma_{1})$ in the obtained punctured disc $D_{n-1}$. That is, we consider a new cutting sequence presentation defined by the intersection with the curve diagram $\Sigma^{(n-1)}$. The second principal part of $\beta$ is, by definition, the $1$st principal part of this cutting sequence. 

The $i$-th principal part is defined similarly. Remove $(n-(i-1)-a_{i-1}(W))$-th puncture point in $D_{n+2-i}$, and shift a numbering of punctures so that the obtained disc is naturally identified with the $(n-(i-1))$ punctured disc $D_{n+1-i}$. The $i$-th principal part of $\beta$ is defined to be the $1$st principal part of the tight cutting sequence presentation of $\beta(\Gamma_{1})$ in the obtained punctured disc $D_{n+1-i}$.   
 
  Geometrically, the $i$-th principal part is the first segment of the arc $\beta(\Gamma_{1})$ which is cut by the subdisc $D'_{n+1-i}$.

\begin{exam}
Let $<_{D}$ be the Dehornoy ordering on $B_{5}$ and consider a dual positive braid $\beta = a_{1,3}a_{3,5}a_{1,5}$.
Then the $\C$-normal form $W=W(\beta;<)$ of $\beta$ is given by
\[
 W  =  A_{1}  A_{0}  A_{-1} = (a_{1,3}a_{3,5}a_{1,5})(\varepsilon )( \varepsilon)
 \]
 and the decomposition of $A_{1}$ is given by 
\[ 
\begin{array}{llc|c|c}
 A_{1}  =  A^{(1)}_{1}  A^{(1)}_{0}  A^{(1)}_{-1} =( a_{1,3}a_{3,5} )( a_{1,5})(\varepsilon).\\
\end{array}
 \] 
 Thus, the $1$st address of $W$ is $1$, and the $2$nd address of $W$ is $1$. 
 
The image $W(\Gamma_{1})$ of the arc $\Gamma_{1}$ is shown in figure \ref{fig:prinexam}.
The $1$st principal part is depicted by the bold line in figure \ref{fig:prinexam} (a).
The $2$nd principal part is depicted by the bold line in figure \ref{fig:prinexam} (b). In figure \ref{fig:prinexam} (b), since $(a_{1}(W))^{*} = 4$, the $4$-th puncture point $p_{4}$ is removed.
Similarly, the $3$rd principal part is depicted by the bold line in the figure \ref{fig:prinexam} (c). In this case, the $2$nd puncture point $p_{2}$ is removed.
 
\begin{figure}[htbp]
 \label{fig:prinexam}
 \begin{center}
\includegraphics[width=100mm]{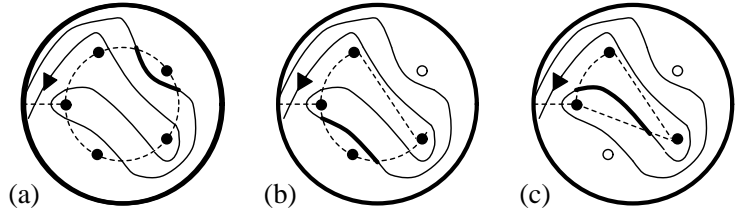}
 \end{center}
 \caption{Principal parts of the braid $a_{1,3}a_{3,5}a_{1,5}$ }

\end{figure}

\end{exam}

\subsection{Proof of theorem \ref{thm:main} and theorem \ref{thm:ordertype}}

Now we prove theorem \ref{thm:main} and \ref{thm:ordertype}. From the definition of the $\C$-normal form for general finite Thurston type orderings, it is sufficient to prove the theorem for normal ordering case. 

Let $<$ be a normal finite Thurston type ordering of $B_{n}$, defined by the normal curve diagram $\Gamma$.
We denote the associated permutation $\{k(1),k(2),\cdots,k(n-1)\}$ of $\Gamma$ by $\bk$.
Throughout this section, we denote $k(1)$ simply by $k$.

First of all, we prove the push-off formula for the dual braid monoids. 

\begin{lem}[push-off formula]
\label{lem:pushoff}
Let $\beta \in B_{\{i,i+1,\cdots,j\}}^{+*} $ and $i\leq j\leq p \leq q$.
then $a_{p,q}\beta = \beta'a_{p,q}$ for some $\beta' \in B_{\{i,i+1,\cdots,q\}}^{+*}$.
 
\end{lem}
\begin{proof}
If $p\neq j$, then $a_{p,q}$ commute with $\beta$, so the assertion is obvious. Assume that $p=j$. Let us write $W= V_{0}a_{i_{1},j}V_{1}\cdots a_{i_{l},j}V_{l}$, where $V_{a} \in B_{\{i,i+1,\cdots,j-1\}}^{+*}$. Then, we obtain 
\[ a_{j,q}\beta = a_{j,q}(V_{0}a_{i_{1},j}V_{1}\cdots a_{i_{l},j}V_{l}) = (V_{0}a_{i_{1},q}V_{1}\cdots a_{i_{l},q}V_{l})a_{j,q}.\]
\end{proof}

Now we proceed to study properties of the $\C$-normal forms.

\begin{lem}
\label{lem:nonempty}
Let $W = A_{m}A_{m-1}\cdots A_{0}A_{-1}$
be a $\C$-normal form and assume that $a_{1}(W)=m$.
Then each word $A_{i}$ is non-empty for $1 \leq i \leq m$.
\end{lem}

\begin{proof}
Since $a_{1}(W)= m$, $A_{m}$ is non-empty word.
Assume $A_{i}$ is an empty word $(1 \leq i < m )$. Then the last letter of $A_{i+1}$ must be $a_{i^{*},i^{*}+1}$ because otherwise, by regarding the last letter of $A_{i+1}$ as the first letter of $A_{i-1}$, we obtain a new decomposition of the word $W$ which has the bigger code.

Since $A_{i-1}\in B_{(i^{*}+1)}^{+*}$, so from lemma \ref{lem:pushoff}, 
$a_{i^{*},i^{*}+1} A_{i-1} = A' a_{ i^{*},i^{*}+1}$ for some dual positive word $A'$. Thus if we regard $a_{i^{*},i^{*}+1}$ as the initial letter of $A_{i-2}$, then we get a new word representative of $W$ which has the bigger code, so it is a contradiction.
\end{proof}

We remark that both $A_{0}$ and $A_{-1}$ might be empty.
For example, the word $W=a_{1,4}||$ is the $\C$-normal form with respect to the Dehornoy ordering $<_{D}$ on $B_{4}$, and both $A_{0}$ and $A_{-1}$ are empty words. 
It is easily observed that such a phenomenon occurs only if $k=1$.

From lemma \ref{lem:nonempty}, each $\C$-normal form can be written as
\[ W = A'_{m}a_{m^{*}+1,j_{m}}| A'_{m-1}a_{(m-1)^{*}+1,j_{m-1}} | \cdots A'_{1}a_{n,j_{1}}| A_{0} |A_{-1}.\]

Now we study a geometric property of $\C$-normal forms.
The following proposition plays a crucial role to prove theorem \ref{thm:main}.  
 
\begin{prop}
\label{prop:keyprop}
Let $W=a_{i,j}^{q}W'$ be a $\C$-normal form and assume $a_{1}(W) \neq -1$.
Then the $d(W)$-th principal part of $W'$ intersects the arc $e_{i,j}$. 
\end{prop}

Before proceeding to the proof of the proposition, we explain the geometrical meanings and consequences of the proposition.

Recall that the $i$-th principal part of a $\C$-normal form $W$ is defined as the first segment of the arc $W(\Gamma_{1})$ which is cut by the subdisc $D'_{n+i-1}$.
Thus roughly speaking, the geometrical meaning of this proposition is that, the arc $W'(\Gamma_{1})$ intersects $e_{i,j}$ at the nearest point from the initial point of $\Gamma_{1}$. In other word, the arc $W'(\Gamma_{1})$ reaches $e_{i,j}$ without any roundabout. Thus, the braid $a_{i,j}$ twists the arc $W'(\gamma_{1})$ at the nearest point.

This proposition allows us to compute the principal parts of the image and provide geometric information about the $\C$-normal forms.
We prepare some new notions. From the definition of the principal part, the $1$st principal part separates the disc $D'=D'_{n}$ into two components. We say the $i$-th principal part {\it lies the right side} of the $1$st principal part if the $i$-th principal part belongs to the right side component of $D'$ or the $i$-th principal part is a subarc of the $1$st principal part. Moreover, if the $i$th principal part is not a subarc of the $1$st principal part, we say the $i$-th principal part {\it strictly lies the right side} of the $1$st principal part.

Now we give consequences of proposition \ref{prop:keyprop}.
\begin{cor}
\label{cor:consequence}
Let 
\[ W = A'_{m}a_{m^{*}+1,j_{m}}| A'_{m-1}a_{(m-1)^{*}+1,j_{m-1}} | \cdots |A'_{1} a_{n,j_{1}}| A_{0}| A_{-1}\] be a $\C$-normal form and assume $m =a_{1}(W)\geq 0$. Then 

\begin{enumerate}
\item The $i$-th principal part of $W$ lies in the right side of the $1$st principal part for all $i$.
\item The initial segment of the tight cutting sequence presentation of $W$ is written as 
\[(\underbrace{0,0,\cdots,0}_{s \textrm{times}},m^{*},-x,\cdots)\]
where $s$ is the maximal integer equal or less than $\frac{m}{n}$. 
Hence the $1$st principal part of $W$ is presented by $(m^{*},-x)$ for some $x$.
\item Moreover, if $m=0$, then $A_{0} \in B_{\{1,2,\cdots,x\}}^{+*}$.
\item Similarly, if $m>1$, then $A'_{m} \in B_{\{m^{*}+1,m^{*}+2,\cdots,x\}}^{+*}$.

\end{enumerate}
\end{cor}

\begin{proof}
We prove lemma by induction on $(n,l)=(n,l(W))$ where $l(W)$ denotes the word length of $W$.
The case $n=2$ or $l=1$ is obvious.
Assume that (1)-(4) is true for all $(n',l')$ with $n'<n$ or $n'=n$ and $l'<l$.
Let $W=a_{p,q}W'$ be the $\C$-normal form of a dual positive $n$-braid whose length is $l$. Throughout the proof we denote the $1$st principal part of $W$ by $\tau$.

Assertion (1) is easily confirmed from our inductive hypothesis. First observe that the half Dehn-twist $a_{p,q}$ preserves the relative positions of arcs. That is, if the arc $\alpha$ in the disc $D'$ lies in the right side of $\tau$, then $a_{p,q}(\alpha)$ also lies in the right side of $a_{p,q}(\tau)$.
Since the $k$-th principal part of $W$ is a subarc of the arc $a_{p,q}(\tau_{l})$ for some $k \leq l$, thus by inductive hypothesis, the $i$-th principal part of $W$ lies in the right side of the $1$st principal part of $W$.

Now we prove the assertion (2), (3) and (4). 
 
Let us consider the case that the length of $A_{m}$ is $1$. 
First assume that $m>0$. Then $A_{m}= a_{m^{*}+1,j_{m}}$. From inductive hypothesis, the cutting sequence presentation of $\tau$ is given by $(m^{*}+1,-x)$ for some $x$, and the arc $e_{m^{*}+1,j_{m}}$ intersects $\tau$. Now by performing $a_{m^{*}+1,j_{m}}$, we obtain that $1$st principal part of $W$ is $(+m^{*},-y)$ for some $y$. (See figure \ref{fig:proof1} (a)).

Next assume that $m=0$. Then $A_{m}=a_{p,q}$ and the cutting sequence presentation of $\tau$ is $(n,-k)$. By proposition, $e_{p,q}$ intersects $\tau$, so the inequality $p<k<q$ holds. Now by performing $a_{p,q}$, we conclude that the $1$st principal part of $W$ is $(n,-q)$. (See figure \ref{fig:proof1} (b)).

\begin{figure}[htbp]
 \label{fig:proof1}
 \begin{center}
\includegraphics[width=110mm]{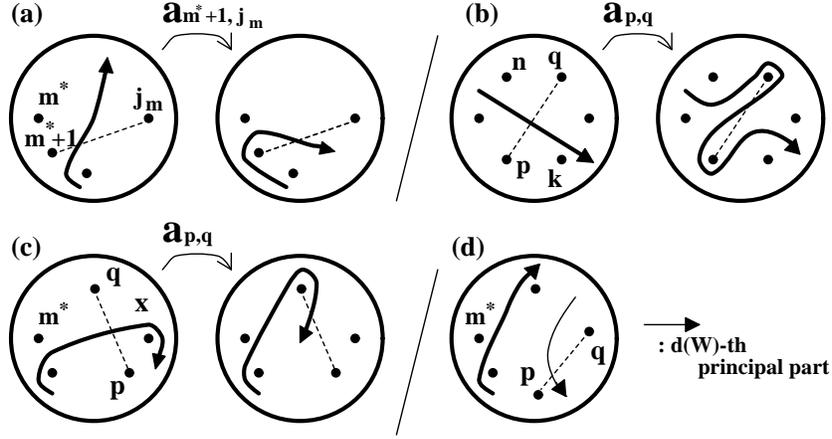}
 \end{center}
 \caption{Change of $\tau$}
\end{figure}

Next we consider the case the length of $A_{m}$ is greater than $1$.
Let us denote the word $A_{m}$ as $A_{m}=a_{p,q}A''a_{m^{*}+1,j_{m}}$. Then from our inductive hypothesis, $\tau = (m^{*},-x)$ for some $x$ and $A'' \in B_{\{m^{*}+1,m^{*}+2,\cdots,x\}}^{+*}$.
If $e_{p,q}$ intersects $\tau$, then $p \leq x < q$ holds. Thus, by performing $a_{p,q}$, we conclude that the $1$st principal part of $W$ is $(m^{*},-q)$. (See figure \ref{fig:proof1} (c)). Thus in this case all assertions (2)-(4) hold.

If $e_{p,q}$ does not intersect $\tau$, then proposition \ref{prop:keyprop} implies $d(W)$ must be greater than $1$. Since $e_{p,q}$ must intersect the $d(W)$-th principal part, which lies in the strictly right side of $\tau$, so the inequality $m^{*} <p<q \leq x$ holds. This implies that $a_{p,q}$ does not change the $1$st principal part $\tau$, so in this case the $1$st principal part of $W$ is also $\tau$ (See figure \ref{fig:proof1} (d)). Thus in this case all assertions (2)-(4) also hold.

\end{proof}

With the above consequences in mind, now we prove proposition \ref{prop:keyprop}.

\begin{proof}[Proof of proposition \ref{prop:keyprop}]
Let $W$ be the $\C$-normal form of a dual positive $n$-braid $\beta$ having the length $l$.
As in the proof of corollary, we prove the proposition by induction on $(n,l)=(n,l(W))$. 
$n=2$ or $l(W)= 0,1$ cases are trivial. 

Now assume the proposition holds for all $(n',l')$ with $n'<n$ or $n=n', l'<l$. We remark that this implies that the results of corollary \ref{cor:consequence} also holds for all $\C$-normal forms for such $(n',l')$.

First we consider the case $d(W)>1$. This case is easily confirmed because it is enough to consider in the sub disc $D_{n-1}$, obtained by filling $a_{1}(W)$-th puncture. From our inductive hypothesis about $n$, we establish the desired result.

Thus, we restrict our attention to the case $d(W)=1$.
Let $m = a_{1}(W)$, the $1$st address of $W$.
If $m=0$, then $W$ is written as $ W= a_{p,q}^{r}A_{-1}$. Since $W$ is the $\C$-normal form, $p \leq k < q$ because otherwise a decomposition of $W$ defined by $A_{-1} = |a_{p,q}^{r}A_{-1}$ has the bigger code. Since the $1$st principal part of $A_{-1}$ is $(+n,-k)$, This implies $e_{p,q}$ intersects the $1$st principal part of $A_{-1}$. 

Thus, from now on, we consider the case $m>0$. 
Denote $W$ as
\[ W=a_{m^{*}+1,j_{m}}^{r}|A'_{m-1}a_{(m-1)^{*}+1,j_{m-1}}| \cdots | A'_{1}a_{n,j_{1}}|A_{0}|A_{-1}\]
and let $\tau$ be the $1$st principal part of $A_{m-1}\cdots A_{-1}$.

Assume that $\tau$ does not intersect $e_{m^{*}+1,j_{m}}$. We deduce a contradiction from this hypothesis.
To this end, we first show the following claim holds under this assumption.

\begin{claim}
Under the above assumption,
\begin{description}
\item[E1(i)] $i^{*}+1<j_{i}<j_{i-1}$
\item[E2(i)] The $1$st principal part of the $\C$-normal form word $A_{i-2}A_{i-3}\cdots A_{-1}$ is $(+i^{*}+1, -j_{i-1})$ 
\end{description}
hold for each $m \geq i \geq 2$.
Moreover, if $A_{0}$ is non-empty, then the $1$st principal part of $A_{-1}$ is $(+n,-j_{0})$ where $A_{0}=A_{0}'a_{n,j_{0}}$.
\end{claim}

\begin{proof}[Proof of claim]

Let us begin with the case $i=m$. 

First observe that from our inductive hypothesis and corollary \ref{cor:consequence}, $\tau = (m^{*}+1,-x)$ for some $x$ and $A'_{m-1} \in B_{\{m^{*}+2,m^{*}+3,\cdots,x\}}^{+*}$.
Now the assumption that $\tau$ does not intersect $e_{m^{*}+1,j_{m}}$ implies that $m^{*}+1<x<j_{m}$ holds. Especially, $j_{m} \neq m^{*}+2$ holds.

If $m^{*}+1<j_{m-1}\leq j_{m}$, then by applying lemma \ref{lem:pushoff} we can push off $a_{m^{*}+1, j_{m}}$ across the word $A'_{m-1}a_{m^{*}+2,j_{m-1}}$. Thus, we can obtain a new word representative 
\[ V= \widetilde{A}'_{m-1} | a_{m^{*}+1,j_{m}}A'_{m-2}a_{m^{*}+3,j_{m-2}}| \cdots \]
where $\widetilde{A}'_{m-1}$ is a dual positive word. Now the code of $V$ is bigger than that of $W$ with respect to the ordering $<_{right}$, which contradicts the assumption $W$ is $\C$-normal form.
Therefore, we conclude that {\bf E1(m)} holds.\\

Now we prove that {\bf E1(i)} and {\bf E2(i+1)} (resp. {\bf E1(m)}) implies {\bf [E2(i)]} (resp. {\bf E2(m)}). Let $\tau_{i}$ be the $1$st principal part of the $\C$-normal form $W=A_{i-2}\cdots A_{-1}$.
From corollary \ref{cor:consequence}, $\tau_{i} = (i^{*}+1,-x')$ for some $x'$.
Let us assume that $x' \neq j_{i-1}$. From our inductive hypothesis, $\tau_{i}$ intersects $e_{i^{*}+2,j_{i-1}}$. Thus the $1$st principal part of $a_{i^{*}+2,j_{i-1}}W$ is given by $(i^{*},-j_{i-1})$. 

On the other hand, from {\bf E2(i+1)} we have already known that the $1$st principal part of $A_{i-1}\cdots A_{-1}$ is given by $(i^{*},-j_{i})$. If $i=m$, then we have also known that $m^{*}+1 < x < j_{m}$. Thus, in either case we obtain the inequality $i^{*} < j_{i-1} \leq j_{i}$. This contradicts {\bf E1(i)}, so we conclude that $\tau_{i}=(i^{*}+1,-j_{i-1})$.\\

Finally we prove that {\bf E2(i)} implies {\bf E1(i-1)}. This complete the proof of claim.
Assume that $(i-1)^{*} +1 = i^{*} +2 < j_{i-2} \leq j_{i-1}$ holds. 
First observe that $j_{i-1} \neq i^{*}+3$ because if so, we obtain $j_{i-2} = i^{*}+2$, which contradicts the fact $A_{i-2} \in B_{(i-2)^{*}}^{+*}$.
From corollary \ref{cor:consequence} and {\bf E2(i)}, $A'_{i-2} \in B_{\{i^{*}+3,i^{*}+4,\cdots,j_{i-1}\}}$ holds. Now by applying lemma \ref{lem:pushoff}, we can push off $a_{i^{*}+2,j_{i-1}}$ across the word $A'_{i-2}a_{i^{*}+3,j_{i-2}}$, and obtain a new word representative 
\[ V' = \cdots \widetilde{A}'_{i-2}| a_{i^{*}+2,j_{i-1}} A'_{i-3} a_{i^{*}+4,j_{i-3}}| \cdots \]
where $\widetilde{A}'_{i-1}$ is a dual positive word. 
This word $V'$ has the bigger code, it is contradiction. Thus we obtain the inequality $m^{*}+2<j_{m-1}<j_{m-2}$.
\end{proof}

Now we deduce a contradiction from the claim.
Suppose that $A_{0}$ is an empty word. Then $k=1$ and $j_{1}=1$ holds.
This is a contradiction because the claim implies $j_{1} \neq 1$.
 If $A_{0}$ is non-empty, then the claim implies that $j_{0}=k$, because the $1$st principal part of $A_{-1}$ is $(+n,-k)$. This contradicts the inductive hypothesis because $e_{p,k}$ does not intersects the $1$st principal part $(+n,-k)$.

This complete the proof for $d(W)=1$ case.

\end{proof}

Now we are ready to prove main theorem.

\begin{proof}[Proof of theorem \ref{thm:main}]

As we mentioned earlier, to prove the theorem, it is sufficient to show the normal finite Thurston type ordering case.
We prove the theorem by induction of $n$. 
The dual $2$-braid monoid $B_{2}^{+*}$ case is obvious.
Assume that we have proved the theorem for normal finite Thurston type orderings on $B_{i}$ for all $i<n$.

 Let $\alpha$ and $\beta$ be dual positive $n$-braids and $<$ be a normal finite Thurston type ordering of $B_{n}$. We denote the $\C$-normal forms of $\alpha$ and $\beta$ with respect to the ordering $<$ by $W$ and $V$ respectively. We assume that $\C(\alpha;<) <_{left} \C(\beta;<)$.

Since Thurston type orderings are left-invariant, we can always assume the initial letter of $W$ and $V$ are different by annihilating their left common divisors. 
If $a_{1}(W)=a_{1}(V)=-1$, then $\alpha<\beta$ follows from the inductive hypothesis. So we assume $a_{1}(V) \geq 0$.

Let $m$ be the minimal positive integer such that the inequality $a_{m}(W) < a_{m}(V)$ holds. Since we assumed that $W$ and $V$ has no left common divisors, we can always find such $m$.

 From the assumption $\C(\alpha;<) <_{left} \C(\beta;<)$, we obtain $a_{i}(V) = a_{i}(W)$ for $i=1,\cdots,m-1$. Thus from corollary \ref{cor:consequence} (2), the $i$-th principal part of $W$ and $V$ are coincide for all $i<m$ and the $i$th principal part of $V$ moves more left than that of $W$.
 Since we can make $W(\Gamma_{1})$ and $V(\Gamma_{1})$ tight without changing the curve diagram $\Sigma$ used to define cutting sequence presentation, this means the arc $V(\Gamma_{1})$ moves more left than $W(\Gamma_{1})$. This implies $\alpha < \beta$ holds.
\end{proof}

Now to compute the order-type is an easy task.

\begin{proof}[Proof of theorem \ref{thm:ordertype}]
For an arbitrary finite Thurston type ordering $<$, let $<_{N}$ be a normal finite Thurston type ordering which is conjugate to $<$ and $\gamma$ be a dual positive conjugating element between $<$ and $<_{N}$.
Then the map $(B_{n}^{+*},<) \rightarrow (B_{n}^{+*};<_{N})$ defined by $\beta \mapsto \beta \gamma$ is order-preserving injection. Thus order-type of $(B_{n}^{+*},<)$ is smaller than that of $(B_{n}^{+*},<_{N})$.
By interchanging the role of $<$ and $<_{N}$, we conclude that $(B_{n}^{+*},<)$ and $(B_{n}^{+*},<_{N})$ have the same order-type. Thus to prove the theorem, it is sufficient to consider the normal ordering case.

Let $<$ be a normal finite Thurston type ordering of $B_{n}$ and $Ncode(n*,<)$ be a set of all codes of the $\C$-normal forms of $B_{n}^{+*}$, with respect to the ordering $<$. Theorem \ref{thm:main} asserts that $(Ncode(n*,<),<_{left})$ is order-isomorphic to $(B_{n}^{+*},<)$.
 A direct computation shows that the order-type of the well-ordered set $(Code(n*,<),<_{left})$ is 
$\omega^{\omega^{n-2}}$. Since $(Ncode(n*,<),<_{left})$ is subset of $(Code(n*,<),<_{left})$, the order type of $(B_{n}^{+*},<)$ is at most $\omega^{\omega^{n-2}}$.

On the other hand, the positive braid monoid $B_{n}^{+}$ is a submonoid of the dual braid monoid $B_{n}^{+*}$, and the order type of $(B_{n}^{+}, <)$ is $\omega^{\omega^{n-2}}$ \cite{i}. Thus the order-type of $(B_{n}^{+*},<)$ is greater than $\omega^{\omega^{n-2}}$. Therefore we conclude that order-type of $(B_{n}^{+*},<)$ is $\omega^{\omega^{n-2}}$.
\end{proof}

We remark that from the above theorem, $(B_{n}^{+},<)$ and $(B_{n}^{+*},<)$ are order-isomorphic.
However, as we remarked in the remark \ref{rem:difference}, the ordinal of a positive braid $\beta$ in the set $(B_{n}^{+*},<)$ might be very different from the ordinal in the set $(B_{n}^{+},<)$.   

\subsection{Property $S$, computational complexity}

Finally we provide applications of theorem \ref{thm:main}.  

First of all we observe that the famous property of Thurston type orderings called the property $S$ (Subword property) can be obtained from the $\C$-normal form description.

\begin{cor}[Property $S$ for finite Thurston type ordering]
Let $<$ be a finite Thurston type ordering on $B_{n}$. Then
\[ \alpha a_{i,j}\beta > \alpha\beta > \alpha a_{i,j}^{-1}\beta \] 
holds for all $\alpha,\beta \in B_{n}$ and $1\leq i<j \leq n-1$.
\end{cor}
\begin{proof}
We only need to prove $a_{i,j} \alpha > \alpha$ for all $\alpha \in B_{n}$. Since $<$ is left-invariant, by multiplying $\delta^{nm}$, which belongs to the center of $B_{n}$ for sufficiently large $m$, we can assume $\alpha \in B_{n}^{+*}$. Let $<'$ be the ordering on $B_{n}$ defined by a conjugate of $<$ by $\alpha^{-1}$. That is, the ordering $<'$ is defined by $\beta >'\gamma$ if and only if $\beta\alpha > \gamma\alpha$ holds. Then $<'$ is also a finite Thurston type ordering. Now from theorem \ref{thm:main}, $a_{i,j} >' 1$ holds so $a_{i,j}\alpha > \alpha$ holds.
\end{proof}

This is interesting in the following sense.
Although we had already established a similar combinatorial description of finite Thurston type orderings using the positive braid monoid $B_{n}^{+}$ in \cite{i}, but in the positive braid monoid $B_{n}^{+}$ case we essentially used the property $S$ to prove the counterpart of theorem \ref{thm:main}. Therefore the $\C$-normal form description for the positive braid monoid does not give a proof of property $S$ although it give a new interpretation of finite Thurston type orderings. 

Finally we prove theorem \ref{thm:complexity}.
The analysis of computational complexity of rotating normal forms and the Dehornoy orderings has been done in \cite{f3}. The computation of a tail twisted rotating normal form is essentially the same as computation of rotating normal form except the tail parts, which are merely computations of the maximal right divisors. 
Since the computations of maximal right divisors are fundamental step of the computation of the rotating normal forms, so the computational complexity of the tail-twisted rotating normal form is identical with the rotating normal forms.

\begin{proof}[Proof of theorem \ref{thm:complexity}]
As we mentioned above, the assertion (1) follows from Fromentin's result \cite{f3}.
We can find two dual positive braid word $V$ and $U$ such that $\beta= U^{-1}V$ in time at most $O(l^{2}n\log n)$ (see \cite{e}, Chapter 9). For the dual braid monoid, the Garside element  $\delta$ has length $n$, so the length of $V$ and $U$ are at most $O(ln)$.
Now $\beta>1$ is equivalent to $V>U$ and the computation of the $\C$-normal form of $V$ and $U$ are at most $O(l^{2}n^{2})$ from the assertion (1). Thus, we conclude that total computational complexity is at most $O(l^{2}n^{2})$.
\end{proof}

 We remark that this computational complexity is smaller than the positive braid monoid case and quadratic with respect to both $n$ and $l$. Moreover, for a given braid word $W$ of length $l$ with respect to Artin generators $\{\sigma_{i}\}$, we can regard $W$ as a braid word of length less than $l$ with respect to the band generators $\{a_{i,j}\}$. Thus, the $\C$-normal form method for dual braid monoid gives more efficient method to compute finite Thurston type orderings.

\end{document}